# On stepdown control of the false discovery proportion

## Joseph P. Romano[1] and Azeem M. Shaikh[2]

*Stanford University*

**Abstract:** Consider the problem of testing multiple null hypotheses. A classical approach to dealing with the multiplicity problem is to restrict attention to procedures that control the familywise error rate ($FWER$), the probability of even one false rejection. However, if $s$ is large, control of the $FWER$ is so stringent that the ability of a procedure which controls the $FWER$ to detect false null hypotheses is limited. Consequently, it is desirable to consider other measures of error control. We will consider methods based on control of the false discovery proportion ($FDP$) defined by the number of false rejections divided by the total number of rejections (defined to be 0 if there are no rejections). The false discovery rate proposed by Benjamini and Hochberg (1995) controls $E(FDP)$. Here, we construct methods such that, for any $\gamma$ and $\alpha$, $P\{FDP > \gamma\} \leq \alpha$. Based on $p$-values of individual tests, we consider stepdown procedures that control the $FDP$, without imposing dependence assumptions on the joint distribution of the $p$-values. A greatly improved version of a method given in Lehmann and Romano [10] is derived and generalized to provide a means by which any sequence of nondecreasing constants can be rescaled to ensure control of the $FDP$. We also provide a stepdown procedure that controls the $FDR$ under a dependence assumption.

## 1. Introduction

In this article, we consider the problem of simultaneously testing a finite number of null hypotheses $H_i$ ($i = 1, \ldots, s$). We shall assume that tests based on $p$-values $\hat{p}_1, \ldots, \hat{p}_s$ are available for the individual hypotheses and the problem is how to combine them into a simultaneous test procedure.

A classical approach to dealing with the multiplicity problem is to restrict attention to procedures that control the familywise error rate ($FWER$), which is the probability of one or more false rejections. In addition to error control, one must also consider the ability of a procedure to detect departures from the null hypotheses when they do occur. When the number of tests $s$ is large, control of the $FWER$ is so stringent that individual departures from the hypothesis have little chance of being detected. Consequently, alternative measures of error control have been considered which control false rejections less severely and therefore provide better ability to detect false null hypotheses.

Hommel and Hoffman [8] and Lehmann and Romano [10] considered the $k$-$FWER$, the probability of rejecting at least $k$ true null hypotheses. Such an error rate with $k > 1$ is appropriate when one is willing to tolerate one or more false rejections, provided the number of false rejections is controlled. They derived single step and stepdown methods that guarantee that the $k$-$FWER$ is bounded above by $\alpha$. Evidently, taking $k = 1$ reduces to the usual $FWER$. Lehmann and Romano [10] also considered control of the *false discovery proportion* ($FDP$), defined as the total number of false rejections divided by the total number of rejections (and equal to 0 if there are no rejections). Given a user specified value $\gamma \in (0, 1)$, control of the $FDP$ means we wish to ensure that $P\{FDP > \gamma\}$ is bounded above by $\alpha$. Control of the *false discovery rate* ($FDR$) demands that $E(FDP)$ is bounded above by $\alpha$. Setting $\gamma = 0$ reduces to the usual $FWER$.

Recently, many methods have been proposed which control error rates that are less stringent than the $FWER$. For example, Genovese and Wasserman [4] study asymptotic procedures that control the $FDP$ (and the $FDR$) in the framework of a random effects mixture model. These ideas are extended in Perone Pacifico, Genovese, Verdinelli and Wasserman [11], where in the context of random fields, the number of null hypotheses is uncountable. Korn, Troendle, McShane and Simon [9] provide methods that control both the $k$-$FWER$ and $FDP$; they provide some justification for their methods, but they are limited to a multivariate permutation model. Alternative methods of control of the $k$-$FWER$ and $FDP$ are given in van der Laan, Dudoit and Pollard [17]. The methods proposed in Lehmann and Romano [10] are not asymptotic and hold under either mild or no assumptions, as long as $p$-values are available for testing each individual hypothesis. In this article, we offer an improved method that controls the $FDP$ under no dependence assumptions of the $p$-values. The method is seen to be a considerable improvement in that the critical values of the new procedure can be increased by typically 50 percent over the earlier procedure, while still maintaining control of the $FDP$. The argument used to establish the improvement is then generalized to provide a means by which any nondecreasing sequence of constants can be rescaled (by a factor that depends on $s$, $\gamma$, and $\alpha$) so as to ensure control of the $FDP$.

It is of interest to compare control of the $FDP$ with control of the $FDR$, and some obvious connections between methods that control the $FDP$ in the sense that

$$P\{FDP > \gamma\} \leq \alpha$$

and methods that control its expected value, the $FDR$, can be made. Indeed, for any random variable $X$ on $[0, 1]$, we have

$$E(X) = E(X|X \leq \gamma)P\{X \leq \gamma\} + E(X|X > \gamma)P\{X > \gamma\}$$
$$\leq \gamma P\{X \leq \gamma\} + P\{X > \gamma\} \,,$$

which leads to

$$(1.1) \qquad \frac{E(X) - \gamma}{1 - \gamma} \leq P\{X > \gamma\} \leq \frac{E(X)}{\gamma},$$

with the last inequality just Markov's inequality. Applying this to $X = FDP$, we see that, if a method controls the $FDR$ at level $q$, then it controls the $FDP$ in the sense $P\{FDP > \gamma\} \leq q/\gamma$. Obviously, this is very crude because if $q$ and $\gamma$ are both small, the ratio can be quite large. The first inequality in (1.1) says that if the $FDP$ is controlled in the sense of (3.3), then the $FDR$ is controlled at level $\alpha(1 - \gamma) + \gamma$, which is $\geq \alpha$ but typically only slightly. Therefore, in principle, a method that controls the $FDP$ in the sense of (3.3) can be used to control the $FDR$ and vice versa.

The paper is organized as follows. In Section 2, we describe our terminology and the general class of stepdown procedures that are examined. Results from Lehmann and Romano [10] are summarized to motivate our choice of critical values. Control of the $FDP$ is then considered in Section 3. The main result is presented in Theorem 3.4 and generalized in Theorem 3.5. In Section 4, we prove that a certain stepdown procedure controls the $FDR$ under a dependence assumption.

## 2. A class of stepdown procedures

A formal description of our setup is as follows. Suppose data $X$ is available from some model $P \in \Omega$. A general hypothesis $H$ can be viewed as a subset $\omega$ of $\Omega$. For testing $H_i : P \in \omega_i$, $i = 1, \ldots, s$, let $I(P)$ denote the set of true null hypotheses when $P$ is the true probability distribution; that is, $i \in I(P)$ if and only if $P \in \omega_i$.

We assume that $p$-values $\hat{p}_1, \ldots, \hat{p}_s$ are available for testing $H_1, \ldots, H_s$. Specifically, we mean that $\hat{p}_i$ must satisfy

$$(2.1) \qquad P\{\hat{p}_i \leq u\} \leq u \quad \text{for any } u \in (0,1) \quad \text{and any } P \in \omega_i,$$

Note that we do not require $\hat{p}_i$ to be uniformly distributed on $(0,1)$ if $H_i$ is true, in order to accomodate discrete situations.

In general, a $p$-value $\hat{p}_i$ will satisfy (2.1) if it is obtained from a nested set of rejection regions. In other words, suppose $S_i(\alpha)$ is a rejection region for testing $H_i$; that is,

$$(2.2) \qquad P\{X \in S_i(\alpha)\} \leq \alpha \quad \text{for all } 0 < \alpha < 1, \ P \in \omega_i$$

and

$$(2.3) \qquad S_i(\alpha) \subset S_i(\alpha') \quad \text{whenever } \alpha < \alpha'.$$

Then, the $p$-value $\hat{p}_i$ defined by

$$(2.4) \qquad \hat{p}_i = \hat{p}_i(X) = \inf\{\alpha : \ X \in S_i(\alpha)\}.$$

satisfies (2.1).

In this article, we will consider the following class of *stepdown* procedures. Let

$$(2.5) \qquad \alpha_1 \leq \alpha_2 \leq \cdots \leq \alpha_s$$

be constants, and let $\hat{p}_{(1)} \leq \cdots \leq \hat{p}_{(s)}$ denote the ordered $p$-values. If $\hat{p}_{(1)} > \alpha_1$, reject no null hypotheses. Otherwise,

$$(2.6) \qquad \hat{p}_{(1)} \leq \alpha_1, \ldots, \hat{p}_{(r)} \leq \alpha_r,$$

and hypotheses $H_{(1)}, \ldots, H_{(r)}$ are rejected, where the largest $r$ satisfying (2.6) is used. That is, a stepdown procedure starts with the most significant $p$-value and continues rejecting hypotheses as long as their corresponding $p$-values are small.

The Holm [6] procedure uses $\alpha_i = \alpha/(s-i+1)$ and controls the $FWER$ at level $\alpha$ under no assumptions on the joint distribution of the $p$-values. Lehmann and Romano [10] generalized the Holm procedure to control the $k$-$FWER$. Specifically, consider the stepdown procedure described in (2.6), where we now take

$$(2.7) \qquad \alpha_i = \begin{cases} \frac{k\alpha}{s} & i \leq k \\ \frac{k\alpha}{s+k-i} & i > k \end{cases}$$

Of course, the $\alpha_i$ depend on $s$ and $k$, but we suppress this dependence in the notation.

**Theorem 2.1 (Hommel and Hoffman [8] and Lehmann and Romano [10]).**
*For testing* $H_i : P \in \omega_i$, $i = 1, \ldots, s$, *suppose* $\hat{p}_i$ *satisfies* (2.1). *The stepdown procedure described in* (2.6) *with* $\alpha_i$ *given by* (2.7) *controls the* $k$-*FWER; that is,*

(2.8)   $P\{\text{reject at least } k \text{ hypotheses } H_i \text{ with } i \in I(P)\} \leq \alpha \quad \text{for all } P .$

*Moreover, one cannot increase even one of the constants* $\alpha_i$ *(for* $i \geq k$*) without violating control of the* $k$-*FWER. Specifically, for* $i \geq k$*, there exists a joint distribution of the* $p$-*values for which*

(2.9)   $P\{\hat{p}_{(1)} \leq \alpha_1, \hat{p}_{(2)} \leq \alpha_2, \ldots, \hat{p}_{(i-1)} \leq \alpha_{i-1}, \hat{p}_{(i)} \leq \alpha_i\} = \alpha.$

**Remark 2.1.** Evidently, one can always reject the hypotheses corresponding to the smallest $k - 1$ $p$-values without violating control of the $k$-FWER. However, it seems counterintuitive to consider a stepdown procedure whose corresponding $\alpha_i$ are not monotone nondecreasing. In addition, automatic rejection of $k - 1$ hypotheses, regardless of the data, appears at the very least a little too optimistic. To ensure monotonicity, our stepdown procedure uses $\alpha_i = k\alpha/s$. Even if we were to adopt the more optimistic strategy of always rejecting the hypotheses corresponding to the $k - 1$ smallest $p$-values, we could still only reject $k$ or more hypotheses if $\hat{p}_{(k)} \leq k\alpha/s$, which is also true for the specific procedure of Theorem 2.1.

## 3. Control of the false discovery proportion

The number $k$ of false rejections that one is willing to tolerate will often increase with the number of hypotheses rejected. So, it might be of interest to control not the number of false rejections (or sometimes called false discoveries) but the proportion of false discoveries. Specifically, let the *false discovery proportion* (*FDP*) be defined by

(3.1)   $FDP = \begin{cases} \frac{\text{Number of false rejections}}{\text{Total number of rejections}} & \text{if the denominator is } > 0 \\ 0 & \text{if there are no rejections} \end{cases}$

Thus *FDP* is the proportion of rejected hypotheses that are rejected erroneously. When none of the hypotheses are rejected, both numerator and denominator of that proportion are 0; since in particular there are no false rejections, the *FDP* is then defined to be 0.

Benjamini and Hochberg [1] proposed to replace control of the *FWER* by control of the *false discovery rate* (*FDR*), defined as

(3.2)   $FDR = E(FDP).$

The *FDR* has gained wide acceptance in both theory and practice, largely because Benjamini and Hochberg proposed a simple stepup procedure to control the *FDR*. Unlike control of the $k$-*FWER*, however, their procedure is not valid without assumptions on the dependence structure of the $p$-values. Their original paper assumed the very strong assumption of independence of $p$-values, but this has been weakened to include certain types of dependence; see Benjamini and Yekutieli [3]. In any case, control of the *FDR* does not prohibit the *FDP* from varying, even if its average value is bounded. Instead, we consider an alternative measure of control that guarantees the *FDP* is bounded, at least with prescribed probability. That is, for a given $\gamma$ and $\alpha$ in $(0, 1)$, we require

(3.3)   $P\{FDP > \gamma\} \leq \alpha.$

To develop a stepdown procedure satisfying (3.3), let $f$ denote the number of false rejections. At step $i$, having rejected $i-1$ hypotheses, we want to guarantee $f/i \leq \gamma$, i.e. $f \leq \lfloor \gamma i \rfloor$, where $\lfloor x \rfloor$ is the greatest integer $\leq x$. So, if $k = \lfloor \gamma i \rfloor + 1$, then $f \geq k$ should have probability no greater than $\alpha$; that is, we must control the number of false rejections to be $\leq k$. Therefore, we use the stepdown constant $\alpha_i$ with this choice of $k$ (which now depends on $i$); that is,

$$(3.4) \qquad \alpha_i = \frac{(\lfloor \gamma i \rfloor + 1)\alpha}{s + \lfloor \gamma i \rfloor + 1 - i}.$$

Lehmann and Romano [10] give two results that show the stepdown procedure with this choice of $\alpha_i$ satisfies (3.3). Unfortunately, some joint dependence assumption on the $p$-values is required. As before, $\hat{p}_1, \ldots, \hat{p}_s$ denotes the $p$-values of the individual tests. Also, let $\hat{q}_1, \ldots, \hat{q}_{|I|}$ denote the $p$-values corresponding to the $|I| = |I(P)|$ true null hypotheses. So $q_i = p_{j_i}$, where $j_1, \ldots, j_{|I|}$ correspond to the indices of the true null hypotheses. Also, let $\hat{r}_1, \ldots, \hat{r}_{s-|I|}$ denote the $p$-values of the false null hypotheses. Consider the following condition: for any $i = 1, \ldots, |I|$,

$$(3.5) \qquad P\{\hat{q}_i \leq u | \hat{r}_1, \ldots, \hat{r}_{s-|I|}\} \leq u;$$

that is, conditional on the observed $p$-values of the false null hypotheses, a $p$-value corresponding to a true null hypothesis is (conditionally) dominated by the uniform distribution, as it is unconditionally in the sense of (2.1). No assumption is made regarding the unconditional (or conditional) dependence structure of the true $p$-values, nor is there made any explicit assumption regarding the joint structure of the $p$-values corresponding to false hypotheses, other than the basic assumption (3.5). So, for example, if the $p$-values corresponding to true null hypotheses are independent of the false ones, but have arbitrary joint dependence within the group of true null hypotheses, the above assumption holds.

**Theorem 3.1 (Lehmann and Romano [10]).** *Assume the condition* (3.5). *Then, the stepdown procedure with $\alpha_i$ given by* (3.4) *controls the FDP in the sense of* (3.3).

Lehmann and Romano [10] also show the same stepdown procedure controls the $FDP$ in the sense of (3.3) under an alternative assumption involving the joint distribution of the $p$-values corresponding to true null hypotheses. We follow their approach here.

**Theorem 3.2 (Lehmann and Romano [10]).** *Consider testing $s$ null hypotheses, with $|I|$ of them true. Let $\hat{q}_{(1)} \leq \cdots \leq \hat{q}_{(|I|)}$ denote the ordered $p$-values for the true hypotheses. Set $M = \min(\lfloor \gamma s \rfloor + 1, |I|)$.*
(i) *For the stepdown procedure with $\alpha_i$ given by* (3.4),

$$(3.6) \qquad P\{FDP > \gamma\} \leq P\{\bigcup_{i=1}^{M}\{\hat{q}_{(i)} \leq \frac{i\alpha}{|I|}\}\}.$$

(ii) *Therefore, if the joint distribution of the $p$-values of the true null hypotheses satisfy Simes inequality; that is,*

$$P\{\{\hat{q}_{(1)} \leq \frac{\alpha}{|I|}\} \bigcup \{\hat{q}_{(2)} \leq \frac{2\alpha}{|I|}\} \bigcup \cdots \bigcup \{\hat{q}_{(|I|)} \leq \alpha\}\} \leq \alpha,$$

*then $P\{FDP > \gamma\} \leq \alpha$.*

Simes inequality is known to hold for many joint distributions of positively dependent variables. For example, Sarkar and Chang [15] and Sarkar [13] have shown that the Simes inequality holds for the family of distributions which is characterized by the multivariate positive of order two condition, as well as some other important distributions.

However, we will argue that the stepdown procedure with $\alpha_i$ given by (3.4) does *not* control the $FDP$ in general. First, we need to recall Lemma 3.1 of Lehmann and Romano [10], stated next for convenience (since we use it later as well). It is related to Lemma 2.1 of Sarkar [13].

**Lemma 3.1.** *Suppose $\hat{p}_1, \ldots, \hat{p}_t$ are p-values in the sense that $P\{\hat{p}_i \leq u\} \leq u$ for all $i$ and $u$ in $(0,1)$. Let their ordered values be $\hat{p}_{(1)} \leq \cdots \leq \hat{p}_{(t)}$. Let $0 = \beta_0 \leq \beta_1 \leq \beta_2 \leq \cdots \leq \beta_m \leq 1$ for some $m \leq t$.*
(i) *Then,*

$$(3.7) \quad P\{\{\hat{p}_{(1)} \leq \beta_1\} \bigcup \{\hat{p}_{(2)} \leq \beta_2\} \bigcup \cdots \bigcup \{\hat{p}_{(m)} \leq \beta_m\}\} \leq t \sum_{i=1}^{m} (\beta_i - \beta_{i-1})/i.$$

(ii) *As long as the right side of (3.7) is $\leq 1$, the bound is sharp in the sense that there exists a joint distribution for the p-values for which the inequality is an equality.*

The following calculation illustrates the fact that the stepdown procedure with $\alpha_i$ given by (3.4) does *not* control the $FDP$ in general.

**Example 3.1.** Suppose $s = 100$, $\gamma = 0.1$ and $|I| = 90$. Construct a joint distribution of $p$-values as follows. Let $\hat{q}_{(1)} \leq \cdots \leq \hat{q}_{(90)}$ denote the ordered $p$-values corresponding to the true null hypotheses. Suppose these 90 $p$-values have some joint distribution (specified below). Then, we construct the $p$-values corresponding to the 10 false null hypotheses conditional on the 90 $p$-values. First, let 8 of the $p$-values corresponding to false null hypotheses be identically zero (or at least less than $\alpha/100$). If $\hat{q}_{(1)} \leq \alpha/92$, let the 2 remaining $p$-values corresponding to false null hypotheses be identically 1; otherwise, if $\hat{q}_{(1)} > \alpha/92$, let the 2 remaining $p$-values also be equal to zero. For this construction, $FDP > \gamma$ if $\hat{q}_{(1)} \leq \alpha/92$ or $\hat{q}_{(2)} \leq 2\alpha/91$. The value of

$$P\{\hat{q}_{(1)} \leq \frac{\alpha}{92} \bigcup \hat{q}_{(2)} \leq \frac{2\alpha}{91}\}$$

can be bounded by Lemma 3.1. The lemma bounds this expression by

$$90 \left( \frac{\alpha}{92} + \frac{\frac{2\alpha}{91} - \frac{\alpha}{92}}{2} \right) \approx 1.48\alpha > \alpha.$$

Moreover, Lemma 3.1 gives a joint distribution for the 90 $p$-values corresponding to true null hypotheses for which this calculation is an equality.

Since one may not wish to assume any dependence conditions on the $p$-values, Lehmann and Romano [10] use Theorem 3.2 to derive a method that controls the $FDP$ without any dependence assumptions. One simply needs to bound the right hand side of (3.6). In fact, Hommel [7] has shown that

$$P\{\bigcup_{i=1}^{|I|} \{\hat{q}_{(i)} \leq \frac{i\alpha}{|I|}\}\} \leq \alpha \sum_{i=1}^{|I|} \frac{1}{i}.$$

This suggests we replace $\alpha$ by $\alpha(\sum_{i=1}^{|I|}(1/i))^{-1}$. But of course $|I|$ is unknown. So one possibility is to bound $|I|$ by $s$ which then results in replacing $\alpha$ by $\alpha/C_s$, where

$$(3.8) \qquad C_j = \sum_{i=1}^{j} \frac{1}{i}.$$

Clearly, changing $\alpha$ in this way is much too conservative and results in a much less powerful method. However, notice in (3.6) that we really only need to bound the union over $M \leq \lfloor \gamma s \rfloor + 1$ events. This leads to the following result.

**Theorem 3.3 (Lehmann and Romano [10]).** *For testing $H_i : P \in \omega_i$, $i = 1, \ldots, s$, suppose $\hat{p}_i$ satisfies (2.1). Consider the stepdown procedure with constants $\alpha_i' = \alpha_i/C_{\lfloor \gamma s \rfloor + 1}$, where $\alpha_i$ is given by (3.4) and $C_j$ defined by (3.8). Then, $P\{FDP > \gamma\} \leq \alpha$.*

The next goal is to improve upon Theorem 3.3. In the definition of $\alpha_i'$, $\alpha_i$ is divided by $C_{\lfloor \gamma s \rfloor + 1}$. Instead, we will construct a stepdown procedure with constants $\alpha_i'' = \alpha_i/D$, where $D = D(\gamma, \alpha, s)$ is much smaller than $C_{\lfloor \gamma s \rfloor + 1}$. This procedure will also control the $FDP$ but, since the critical values $\alpha_i''$ are uniformly bigger than the $\alpha_i'$, the new procedure can reject more hypotheses and hence is more powerful.

To this end, define

$$(3.9) \qquad \beta_m = \frac{m}{\max\{s + m - \lceil \frac{m}{\gamma} \rceil + 1, |I|\}} \qquad m = 1, \ldots, \lfloor \gamma s \rfloor$$

and

$$(3.10) \qquad \beta_{\lfloor \gamma s \rfloor + 1} = \frac{\lfloor \gamma s \rfloor + 1}{|I|}.$$

where $\lceil x \rceil$ is the least integer $\geq x$.

Next, let

$$(3.11) \qquad N = N(\gamma, s, |I|) = \min\{\lfloor \gamma s \rfloor + 1, |I|, \lfloor \gamma(\frac{s - |I|}{1 - \gamma} + 1) \rfloor + 1\}.$$

Then, let $\beta_0 = 0$ and set

$$(3.12) \qquad S = S(\gamma, s, |I|) = |I| \sum_{i=1}^{N} \frac{\beta_i - \beta_{i-1}}{i}.$$

Finally, let

$$(3.13) \qquad D = D(\gamma, s) = \max_{|I|} S(\gamma, s, |I|).$$

**Theorem 3.4.** *For testing $H_i : P \in \omega_i$, $i = 1, \ldots, s$, suppose $\hat{p}_i$ satisfies (2.1). Consider the stepdown procedure with constants $\alpha_i'' = \alpha_i/D(\gamma, s)$, where $\alpha_i$ is given by (3.4) and $D(\gamma, s)$ is defined by (3.13). Then, $P\{FDP > \gamma\} \leq \alpha$.*

*Proof.* Let $\alpha'' = \alpha/D$. Denote by

$$\hat{q}_{(1)} \leq \cdots \leq \hat{q}_{(|I|)}$$

the ordered $p$-values corresponding only to true null hypotheses. Let $j$ be the smallest (random) index where the $FDP$ exceeds $\gamma$ for the first time at step $j$; that is, the number of false rejections out of the first $j-1$ rejections divided by $j$ exceeds $\gamma$ for the first time at $j$. Denote by $m > 0$ the unique integer satisfying $m - 1 \leq \gamma j < m$. Then, at step $j$, it must be the case that $m$ true null hypotheses have been rejected. Hence,

$$\hat{q}_{(m)} \leq \alpha_j'' = \frac{m\alpha''}{s + m - j}.$$

Note that the number of true hypotheses $|I|$ satisfies

$$|I| \leq s + m - j.$$

Further note that $\gamma j < m$ implies that

(3.14) $$j \leq \lceil \frac{m}{\gamma} \rceil - 1.$$

Hence, $\alpha_j''$ is bounded above by $\beta_m$ defined by (3.9) whenever $m - 1 \leq \gamma j < m$. Note that, when $m = \lfloor \gamma s \rfloor + 1$, we bound $\alpha_j''$ by using $j \leq s$ rather than (3.14).

The possible values of $m$ that must be considered can be bounded. First of all, $j \leq s$ implies that $m \leq \lfloor \gamma s \rfloor + 1$. Likewise, it must be the case that $m \leq |I|$. Finally, note that $j > \frac{s-|I|}{1-\gamma}$ implies that $FDP > \gamma$. To see this, observe that

$$\frac{s - |I|}{1 - \gamma} = (s - |I|) + \frac{\gamma}{1 - \gamma}(s - |I|),$$

so at such a step $j$, it must be the case that

$$t > \frac{\gamma}{1 - \gamma}(s - |I|)$$

true null hypotheses have been rejected. If we denote by $f = j - t$ the number of false null hypotheses that have been rejected at step $j$, it follows that

$$t > \frac{\gamma}{1 - \gamma}f,$$

which in turn implies that

$$FDP = \frac{t}{t + f} > \gamma.$$

Hence, for $j$ to satisfy the above assumption of minimality, it must be the case that

$$j - 1 \leq \frac{s - |I|}{1 - \gamma},$$

from which it follows that we must also have

$$m \leq \lfloor \gamma(\frac{s - |I|}{1 - \gamma} + 1) \rfloor + 1.$$

Therefore, with $N$ defined in (3.11) and $j$ defined as above, we have that

$$P\{FDP > \gamma\} \leq \sum_{m=1}^{N} P\left\{ \{\hat{q}_{(m)} \leq \alpha_j''\} \bigcap \{m - 1 \leq \gamma j < m\} \right\}$$

$$\leq \sum_{m=1}^{N} P\left\{ \hat{q}_{(m)} \leq \alpha''\beta_m \right\} \bigcap \{m - 1 \leq \gamma j < m\} \right\}$$

$$\leq \sum_{m=1}^{N} P\left\{\bigcup_{i=1}^{N}\{\hat{q}_{(i)} \leq \alpha'' \beta_i\} \bigcap \{m-1 \leq \gamma j < m\}\right\}$$

$$\leq P\left\{\bigcup_{i=1}^{N}\{\hat{q}_{(i)} \leq \alpha'' \beta_i\}\right\}.$$

Note that $\beta_m \leq \beta_{m+1}$. To see this, observed that the expression $m + s - \lceil \frac{m}{\gamma} \rceil + 1$ is monotone nonincreasing in $m$, and so the denominator of $\beta_m$, $\max\{m + s - \lceil \frac{m}{\gamma} \rceil + 1, |I|\}$, is monotone nonincreasing in $m$ as well. Also observe that $\beta_m \leq m/|I| \leq 1$ whenever $m \leq N$. We can therefore apply Lemma 3.1 to conclude that

$$P\{FDP > \gamma\} \leq \alpha'' |I| \sum_{i=1}^{N} \frac{\beta_i - \beta_{i-1}}{i}$$

$$= \frac{\alpha |I|}{D} \sum_{i=1}^{N} \frac{\beta_i - \beta_{i-1}}{i} = \frac{\alpha S}{D} \leq \alpha,$$

where $S$ and $D$ are defined in (3.12) and (3.13), respectively. $\qquad \square$

It is important to note that by construction the quantity $D(\gamma, s)$, which is defined to be the maximum over the possible values of $|I|$ of the quantity $S(\gamma, s, |I|)$, does not depend on the unknown number of true hypotheses. Indeed, if the number of true hypotheses, $|I|$, were known, then the smaller quantity $S(\gamma, s, |I|)$ could be used in place of $D(\gamma, s)$.

Unfortunately, a convenient formula is not available for $D(\gamma, s)$, though it is simple to program its evaluation. For example, if $s = 100$ and $\gamma = 0.1$, then $D = 2.0385$. In contrast, the constant $C_{\lfloor \gamma s \rfloor + 1} = C_{11} = 3.0199$. In this case, the value of $|I|$ that maximizes $S$ to yield $D$ is 55. Below, in Table 1 we evaluate $D(\gamma, s)$ and $C_{\lfloor \gamma s \rfloor + 1}$ for several different values of $\gamma$ and $s$. We also compute the ratio of $C_{\lfloor \gamma s \rfloor + 1}$ to $D(\gamma, s)$, from which it is possible to see the magnitude of the improvement of the Theorem 3.4 over Theorem 3.3: the constants of Theorem 3.4 are generally about 50 percent larger than those of Theorem 3.3.

**Remark 3.1.** The following crude argument suggests that, for critical values of the form $d\alpha_i$ for some constant $d$, the value of $d = D^{-1}(\gamma, s)$ is very nearly the largest possible constant one can use and still maintan control of the $FDP$. Consider the case where $s = 1000$ and $\gamma = .1$. In this instance, the value of $|I|$ that maximizes $S$ is 712, yielding $N = 33$ and $D = 3.4179$. Suppose that $|I| = 712$ and construct the joint distribution of the 288 $p$-values corresponding to false hypotheses as follows: For $1 \leq i \leq 28$, if $\hat{q}_{(i)} \leq \alpha \beta_i$ and $\hat{q}_{(j)} > \alpha \bar{\beta}_j$ for all $j < i$, then let $\lceil \frac{i}{\gamma} \rceil - 1$ of the false $p$-values be 0 and set the remainder equal to 1. Let the joint distribution of the 712 true $p$-values be constructed according to the configuration in Lemma 3.1. Note that for such a joint distribution of $p$-values, we have that

$$P\{FDP > \gamma\} \geq P\left\{\bigcup_{i=1}^{28}\{\hat{q}_i \leq \alpha \beta_i\}\right\} = \alpha |I| \sum_{i=1}^{28} \frac{\beta_i - \beta_{i-1}}{i} = 3.2212\alpha.$$

Hence, the largest one could possibly increase the constants by a multiple and still maintain control of the $FDP$ is by a factor of $3.4179/3.2212 \approx 1.061$.

TABLE 1
*Values of $D(\gamma, s)$ and $C_{\lfloor \gamma s \rfloor + 1}$*

| $s$ | $\gamma$ | $D(\gamma, s)$ | $C_{\lfloor \gamma s \rfloor + 1}$ | Ratio |
|------|------|--------|--------|--------|
| 100 | 0.01 | 1 | 1.5 | 1.5 |
| 250 | 0.01 | 1.4981 | 1.8333 | 1.2238 |
| 500 | 0.01 | 1.7246 | 2.45 | 1.4206 |
| 1000 | 0.01 | 2.0022 | 3.0199 | 1.5083 |
| 2000 | 0.01 | 2.3515 | 3.6454 | 1.5503 |
| 5000 | 0.01 | 2.8929 | 4.5188 | 1.562 |
| 25 | 0.05 | 1.4286 | 1.5 | 1.05 |
| 50 | 0.05 | 1.4952 | 1.8333 | 1.2262 |
| 100 | 0.05 | 1.734 | 2.45 | 1.4129 |
| 250 | 0.05 | 2.1237 | 3.1801 | 1.4974 |
| 500 | 0.05 | 2.4954 | 3.8544 | 1.5446 |
| 1000 | 0.05 | 2.9177 | 4.5188 | 1.5488 |
| 2000 | 0.05 | 3.3817 | 5.1973 | 1.5369 |
| 5000 | 0.05 | 4.0441 | 6.1047 | 1.5095 |
| 10 | 0.1 | 1 | 1.5 | 1.5 |
| 25 | 0.1 | 1.4975 | 1.8333 | 1.2242 |
| 50 | 0.1 | 1.7457 | 2.45 | 1.4034 |
| 100 | 0.1 | 2.0385 | 3.0199 | 1.4814 |
| 250 | 0.1 | 2.5225 | 3.8544 | 1.528 |
| 500 | 0.1 | 2.9502 | 4.5188 | 1.5317 |
| 1000 | 0.1 | 3.4179 | 5.1973 | 1.5206 |
| 2000 | 0.1 | 3.9175 | 5.883 | 1.5017 |
| 5000 | 0.1 | 4.6154 | 6.7948 | 1.4722 |

It is worthwhile to note that the argument used in the proof of Theorem 3.4 does not depend on the specific form of the original $\alpha_i$. In fact, it can be used with any nondecreasing sequence of constants to construct a stepdown procedure that controls the $FDP$ by scaling the constants appropriately. To see that this is the case, consider any nondecreasing sequence of constants $\delta_1 \leq \cdots \leq \delta_s$ such that $0 \leq \delta_i \leq 1$ (this restriction is without loss of generality since it can always be acheived by rescaling the constants if necessary) and redefine the constants $\beta_m$ of equations (3.9) and (3.10) by the rule

$$(3.15) \qquad \beta_m = \delta_{k(s,\gamma,m,|I|)} \qquad m = 1, \ldots, \lfloor \gamma s \rfloor + 1$$

where

$$k(s, \gamma, m, |I|) = \min\{s, s + m - |I|, \lceil \frac{m}{\gamma} \rceil - 1\}.$$

Note that in the special case where $\delta_i = \alpha_i$, the definition of $\beta_m$ in equation (3.15) agrees with the earlier definition of equations (3.9) and (3.10). Maintaining the definitions of $N$, $S$, and $D$ in equations (3.11) - (3.13) (where they are now defined in terms of the $\beta_m$ sequence given by equation (3.15)), we then have the following result:

**Theorem 3.5.** *For testing $H_i : P \in \omega_i$, $i = 1, \ldots, s$, suppose $\hat{p}_i$ satisfies (2.1). Let $\delta_1 \leq \cdots \leq \delta_s$ be any nondecreasing sequence of constants such that $0 \leq \delta_i \leq 1$ and consider the stepdown procedure with constants $\delta_i'' = \alpha \delta_i / D(\gamma, s)$, where $D(\gamma, s)$ is defined by (3.13). Then, $P\{FDP > \gamma\} \leq \alpha$.*

*Proof.* Define $j$ and $m$ as in the proof of Theorem 3.4. We have, as before, that whenever $m - 1 \leq \gamma j < m$

$$|I| \leq s + m - j,$$

and

$$j \leq \lceil \frac{m}{\gamma} \rceil - 1.$$

Since $j \leq s$, it follows that

$$\hat{q}_{(m)} \leq \delta_j \leq \beta_m,$$

where $\beta_m$ is as defined in (3.15). The remainder of the argument is identical to the proof of Theorem 3.4 so we do not repeat it here. □

As an illustration of this more general result, consider the nondecreasing sequence of constants given simply by $\eta_i = \frac{i}{s}$. These constants are proportional to the constants used in the procedures for controlling the $FDR$ by Benjamini and Hochberg [1] and Benjamini and Yekutieli [3]. Applying Theorem 3.5 to this sequence of constants yields the following corollary:

**Corollary 3.1.** *For testing $H_i : P \in \omega_i$, $i = 1, \ldots, s$, suppose $\hat{p}_i$ satisfies (2.1). Then the following are true:*
(i) *The stepdown procedure with constants $\eta_i' = \alpha\eta_i/D(\gamma, s)$, where $D(\gamma, s)$ is defined by (3.13), satisfies $P\{FDP > \gamma\} \leq \alpha$;*
(ii) *The stepdown procedure with constants $\eta_i'' = \gamma\alpha\eta_i/\max\{C_{\lfloor\gamma s\rfloor}, 1\}$, where $C_0$ is understood to equal 0, satisfies $P\{FDP > \gamma\} \leq \alpha$.*

*Proof.* The proof of (i) follows immediately from Theorem 3.5. To prove (ii), first observe that $N \leq \lfloor\gamma s\rfloor + 1$ and that for this particular sequence, we have that $\beta_m \leq \min\{\frac{m}{\gamma s}, 1\} =: \zeta_m$. Hence, we have that

$$P\{\bigcup_{i=1}^{N}\{\hat{q}_{(m)} \leq \beta_m\}\} \leq P\{\bigcup_{m=1}^{\lfloor\gamma s\rfloor+1}\{\hat{q}_{(m)} \leq \zeta_m\}\}.$$

Using Lemma 3.1, we can bound the righthand side of this inequality by the sum

$$|I| \sum_{m=1}^{\lfloor\gamma s\rfloor+1} \frac{\zeta_m - \zeta_{m-1}}{m}.$$

Whenever $\lfloor\gamma s\rfloor \geq 1$, we have that $\zeta_{\lfloor\gamma s\rfloor+1} = \zeta_{\lfloor\gamma s\rfloor} = s$, so this sum can in turn be bounded by

$$\frac{|I|}{\gamma s} \sum_{m=1}^{\lfloor\gamma s\rfloor} \frac{1}{m} \leq \frac{1}{\gamma}C_{\lfloor\gamma s\rfloor}.$$

If, on the other hand, $\lfloor\gamma s\rfloor = 0$, we can simply bound the sum by $\frac{1}{\gamma}$. Therefore, if we let $C_0 = 0$, we have that

$$D(\gamma, s) \leq \frac{1}{\gamma}\max\{C_{\lfloor\gamma s\rfloor}, 1\},$$

from which the desired claim follows. □

In summary, given any nondecreasing sequence of constants $\delta_i$, we have derived a stepdown procedure which controls the $FDP$, and so it is interesting to compare such $FDP$-controlling procedures. Clearly, a procedure with larger critical values is preferable to one with smaller ones, subject to the error constraint. The discussion from Remark 3.1 leads us to believe that the critical values from a single procedure will not uniformly dominate those from another, at least approximately. We now consider some specific comparisons which may shed light on how to choose among the various procedures.

TABLE 2
*Values of $D(\gamma, s)$ and $\frac{1}{\gamma}\max\{C_{\lfloor\gamma s\rfloor}, 1\}$*

| $s$ | $\gamma$ | $D(\gamma, s)$ | $\frac{1}{\gamma}\max\{C_{\lfloor\gamma s\rfloor}, 1\}$ | Ratio |
|-----|----------|----------------|------------------|-------|
| 100 | 0.01 | 25.5 | 100 | 3.9216 |
| 250 | 0.01 | 60.4 | 150 | 2.4834 |
| 500 | 0.01 | 90.399 | 228.33 | 2.5258 |
| 1000 | 0.01 | 128.53 | 292.9 | 2.2788 |
| 2000 | 0.01 | 171.73 | 359.77 | 2.095 |
| 5000 | 0.01 | 235.94 | 449.92 | 1.9069 |
| 25 | 0.05 | 6.76 | 20 | 2.9586 |
| 50 | 0.05 | 12.4 | 30 | 2.4194 |
| 100 | 0.05 | 18.393 | 45,667 | 2.4828 |
| 250 | 0.05 | 28.582 | 62.064 | 2.1714 |
| 500 | 0.05 | 37.513 | 76.319 | 2.0345 |
| 1000 | 0.05 | 47.26 | 89.984 | 1.904 |
| 2000 | 0.05 | 57.666 | 103.75 | 1.7991 |
| 5000 | 0.05 | 72.126 | 122.01 | 1.6917 |
| 10 | 0.1 | 3 | 10 | 3.3333 |
| 25 | 0.1 | 6.4 | 15 | 2.3438 |
| 50 | 0.1 | 9.3867 | 22.833 | 2.4325 |
| 100 | 0.1 | 13.02 | 29.29 | 2.2496 |
| 250 | 0.1 | 18.834 | 38.16 | 2.0261 |
| 500 | 0.1 | 23.703 | 44.992 | 1.8981 |
| 1000 | 0.1 | 28.886 | 51.874 | 1.7958 |
| 2000 | 0.1 | 34.317 | 58.78 | 1.7129 |
| 5000 | 0.1 | 41.775 | 67.928 | 1.6261 |

To compare the constants from parts (i) and (ii) of Corollary 3.1, Table 2 displays $D(\gamma, s)$ and $\frac{1}{\gamma}\max\{C_{\lfloor\gamma s\rfloor}, 1\}$ for several different values of $s$ and $\gamma$, as well as the ratio $\frac{1}{\gamma}\max\{C_{\lfloor\gamma s\rfloor}, 1\}/D(\gamma, s)$. In this instance, the improvement between the constants from part (i) and part (ii) is dramatic: The constants $\eta_i'$ are often at least twice as large as the constants $\eta_i''$.

It is also of interest to compare the constants from part (i) of the corollary with those from Theorem 3.4. We do this for the case in which $s = 100$, $\gamma = .1$, and $\alpha = .05$ in Figure 1. The top panel displays the constants $\alpha_i''$ from Theorem 3.4 and the middle panel displays the constants $\eta_i'$ from Corollary 3.1 (i). Note that the scale of the top panel is much larger than the scale of the middle panel. It is therefore clear that the constants $\alpha_i''$ are generally much larger than the constants $\eta_i'$. But it is important to note that the constants from Theorem 3.4 are not uniformly larger than the constants from Corollary 3.1 (i). To make this clear, the bottom panel of Figure 1 displays the ratio $\alpha_i''/\eta_i'$. Notice that at steps 7 - 9, 15 - 19, and 25 - 29 the ratios are strictly less than 1, meaning that at those steps the $\eta_i'$ are larger than the $\alpha_i''$. Following our discussion in Remark 3.1 that these constants are nearly the best possible up to a scalar multiple, we should expect that this would be the case because otherwise the constants $\eta_i'$ could be multiplied by a factor larger than 1 and still retain control of the *FDP*. Even at these steps, however, the constants $\eta_i'$ are very close to the constants $\alpha_i''$ in absolute terms. Since the constants $\alpha_i''$ are considerably larger than the constants $\eta_i'$ at other steps, this suggests that the procedure based upon the constants $\alpha_i''$ is preferable to the procedure based on the constants $\eta_i'$.

## 4. Control of the *FDR*

Next, we construct a stepdown procedure that controls the FDR under the same conditions as Theorem 3.1. The dependence condition used is much weaker than that of independence of $p$-values used by Benjamini and Liu [2].

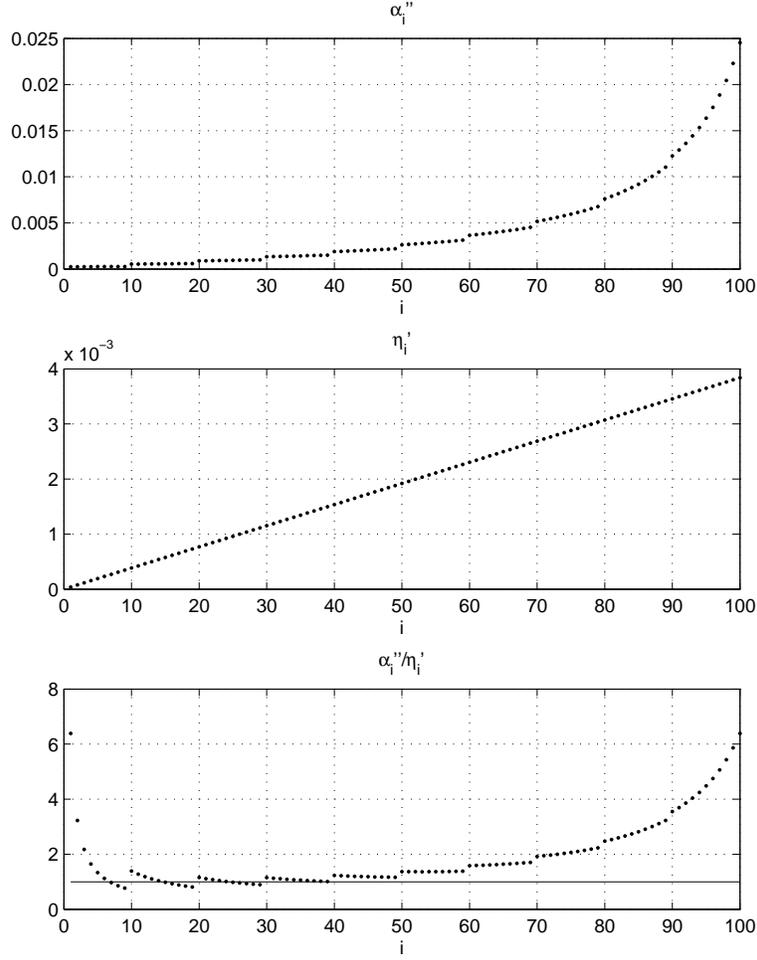

FIG 1. *Stepdown Constants for* $s = 100$, $\gamma = .1$, *and* $\alpha = .05$.

**Theorem 4.1.** *For testing* $H_i : P \in \omega_i$, $i = 1, \ldots, s$, *suppose* $\hat{p}_i$ *satisfies* (2.1). *Consider the stepdown procedure with constants*

$$(4.1) \qquad \alpha_i^* = \min\{\frac{s\alpha}{(s-i+1)^2}, 1\}$$

*and assume the condition* (3.5). *Then,* $FDR \leq \alpha$.

*Proof.* First note that if $|I| = 0$, then $FDR = 0$. Second, if $|I| = s$, then $FDR = P\{\hat{p}_{(1)} \leq \alpha_1^*\} \leq \sum_{i=1}^s P\{\hat{p}_i \leq \alpha_1^*\} \leq s\alpha_1^* = \alpha$.

Now suppose that $0 < |I| < s$. Define $\hat{q}_1, \ldots, \hat{q}_{|I|}$ and $\hat{r}_1, \ldots, \hat{r}_{s-|I|}$ to be the $p$-values corresponding, respectively, to the true and false hypotheses, and let $\hat{q}_{(1)} \leq \cdots \leq \hat{q}_{(|I|)}$ and $\hat{r}_{(1)} \leq \cdots \leq \hat{r}_{(s-|I|)}$ be their ordered values. Denote by $j$ the largest index such that $\hat{r}_{(1)} \leq \alpha_1^*, \ldots, \hat{r}_{(j)} \leq \alpha_j^*$ (defined to be 0 if $\hat{r}_{(1)} > \alpha_1^*$). Define $t$ to be the total number of true hypotheses rejected by the stepdown procedure and $f$ to be the total number of false hypotheses rejected by the stepdown procedure.

Using this notation, observe that

$$E(FDP|\hat{r}_1, \ldots, \hat{r}_{s-|I|}) = E(\frac{t}{t+f}\{t+f>0\}|\hat{r}_1, \ldots, \hat{r}_{s-|I|})$$

$$\leq E(\frac{t}{t+j}\{t>0\}|\hat{r}_1, \ldots, \hat{r}_{s-|I|})$$

$$\leq \frac{|I|}{|I|+j}E(\{t>0\}|\hat{r}_1, \ldots, \hat{r}_{s-|I|})$$

$$\leq \frac{|I|}{|I|+j}P\{\hat{q}_{(1)} \leq \alpha^*_{j+1}|\hat{r}_1, \ldots, \hat{r}_{s-|I|}\}$$

$$\leq \frac{|I|}{|I|+j}\sum_{i=1}^{|I|}P\{\hat{q}_i \leq \alpha^*_{j+1}|\hat{r}_1, \ldots, \hat{r}_{s-|I|}\}$$

(4.2)
$$\leq \frac{|I|}{|I|+j}|I|\alpha^*_{j+1}$$

$$\leq \frac{|I|^2}{|I|+j}\min\{\frac{s\alpha}{(s-j)^2}, 1\}$$

(4.3)
$$\leq \frac{|I|\alpha}{(s-j)}\frac{|I|s}{(|I|+j)(s-j)}.$$

The inequality (4.2) follows from the assumption (3.5) on the joint distribution of $p$-values. To complete the proof, note that $|I|+j \leq s$. It follows that $\frac{|I|\alpha}{(s-j)} \leq \alpha$ and $(|I|+j)(s-j) - |I|s = j(s-|I|) - j^2 = j(s-|I|-j) \geq 0$. Combining these two inequalities, we have that the expression in (4.3) is bounded above by $\alpha$. The desired bound for the $FDR$ follows immediately. $\qquad\square$

The following simple example illustrates the fact that the $FDR$ is not controlled by the stepdown procedure with constants $\alpha^*_i$ absent the restriction (3.5) on the dependence structure of the $p$-values.

**Example 4.1.** Suppose there are $s = 3$ hypotheses, two of which are true. In this case, $\alpha^*_1 = \frac{\alpha}{3}$, $\alpha^*_2 = \frac{3\alpha}{4}$, and $\alpha^*_3 = \min\{3\alpha, 1\}$. Define the joint distribution of the two true $p$-values $q_1$ and $q_2$ as follows: Denote by $I_i$ the half open interval $[\frac{i-1}{3}, \frac{i}{3})$ and let $(q_1, q_2) \sim U(I_i \times I_j)$ with probability $\frac{1}{6}$ for all $(i,j)$ such that $i \neq j$, $1 \leq i \leq 3$ and $1 \leq j \leq 3$. It is easy to see that $(q_{(1)}, q_{(2)}) \sim U(I_i \times I_j)$ with probability $\frac{1}{3}$ for all $(i,j)$ such that $i < j$, $1 \leq i \leq 3$ and $1 \leq j \leq 3$. Now define the distribution of the false $p$-value $r_1$ conditional on $(q_1, q_2)$ by the following rule: If $q_{(1)} \leq \alpha/3$, then let $r_1 = 1$; otherwise, let $r_1 = 0$. For such a joint distribution of $(q_1, q_2, r_1)$, we have that the $FDP$ is identically one whenever $q_{(1)} \leq \frac{\alpha}{3}$ and is at least $\frac{1}{2}$ whenever $\frac{\alpha}{3} < q_{(1)} \leq \frac{3\alpha}{4}$. Hence,

$$FDR \geq P\{q_{(1)} \leq \frac{\alpha}{3}\} + \frac{1}{2}P\{\frac{\alpha}{3} < q_{(1)} \leq \frac{3\alpha}{4}\}.$$

For $\alpha < \frac{4}{9}$, we therefore have that

$$FDR \geq \frac{2\alpha}{3} + (\frac{3\alpha}{4} - \frac{\alpha}{3}) = \frac{13\alpha}{12} > \alpha.$$

**Remark 4.1.** Some may find it unpalatable to allow the constants to exceed $\alpha$. In this case, one might consider replacing the constants $\alpha_i^*$ above with the more conservative values $\alpha \min\{\frac{s}{(s-i+1)^2}, 1\}$, which by construction are always less than $\alpha$. Since these constants are uniformly smaller than the $\alpha_i^*$, our method of proof shows that the $FDR$ would still be controlled under the dependence condition (3.5). The above counterexample, which did not depend on the particular value of $\alpha_3^*$, however, would show that it is not controlled in general.

Under the dependence condition (3.5), the constants (4.1) control the $FDR$ in the sense $FDR \leq \alpha$, while the constants given by (3.4) control the $FDP$ in the sense of (3.3). Utilizing (1.1), we can use the constants (4.1) to control the $FDP$ by controlling the $FDR$ at level $\alpha\gamma$. In Figure 2, we plot the constants (3.4) and (4.1) for the special case in which $s = 100$ and we use both constants to control the $FDP$ for $\gamma = .1$, and $\alpha = .05$.

The top panel displays the constants $\alpha_i$, the middle panel displays the constants $\alpha_i^*$, and the bottom panel displays the ratio $\alpha_i/\alpha_i^*$. Since the ratios essentially always exceed 1, it is clear that in this instance the constants (3.4) are superior to the constants (4.1). If by utilizing (1.1) we use the constants (3.4) to control the $FDR$, on the other hand, we find that the reverse is true. Control of the $FDR$ at level $\alpha$ can be achieved, for example, by controlling the $FDP$ at level $\frac{\alpha}{2-\alpha}$ and letting $\gamma = \frac{\alpha}{2}$. Figure 3 plots the constants (3.4) and (4.1) for the special case in which $s = 100$ and we use both constants to control the $FDR$ at level $\alpha = .05$.

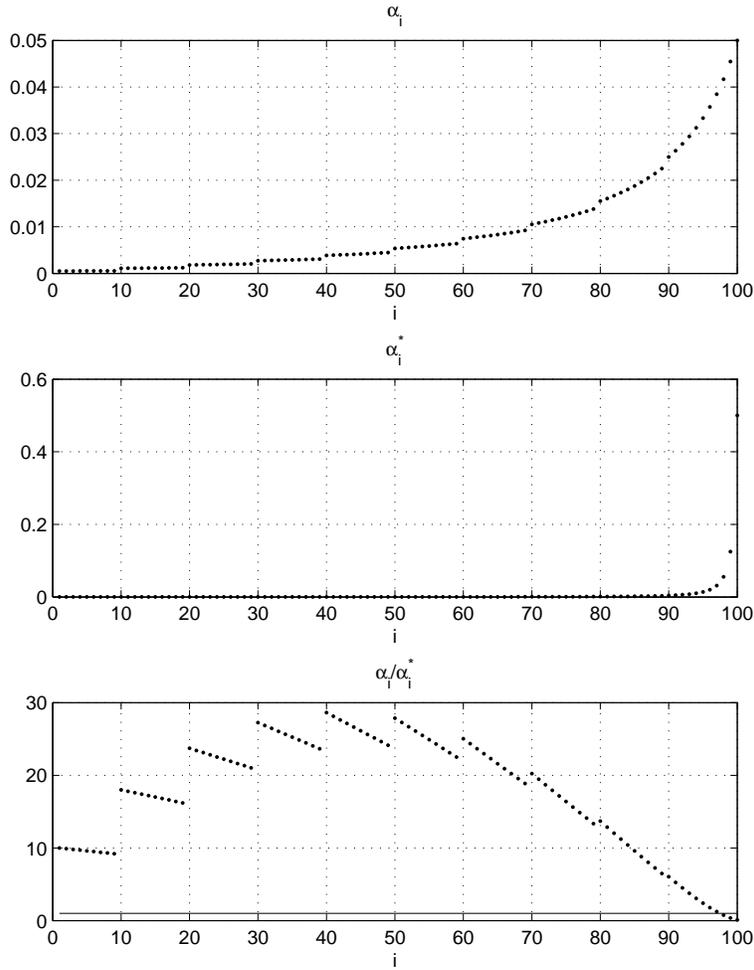

Fig 2. *FDP Control for $s = 100$, $\gamma = .1$, and $\alpha = .05$.*

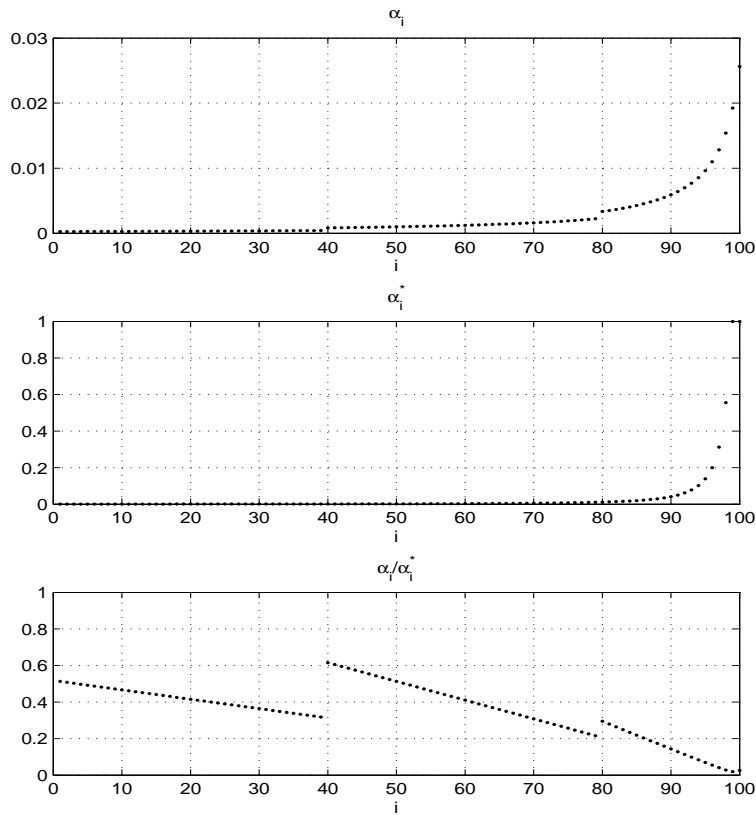

FIG 3. *FDR Control for $s = 100$ and $\alpha = .05$.*

As before, the top panel displays the constants $\alpha_i$, the middle panel displays the constants $\alpha_i^*$, and the bottom panel displays the ratio $\alpha_i/\alpha_i^*$. In this case, the ratio is always less than 1. Thus, in this instance, the constants $\alpha_i^*$ are preferred to the constants $\alpha_i$. Of course, the argument used to establish (1.1) is rather crude, but it nevertheless suggests that it is worthwhile to consider the type of control desired when choosing critical values.

## 5. Conclusions

In this article we have described stepdown procedures for testing multiple hypotheses that control the $FDP$ without any restrictions on the joint distribution of the $p$-values. First, we have improved upon a method proposed by Lehmann and Romano [10]. The new procedure is a considerable improvement in the sense that its critical values are generally 50 percent larger than those of the earlier procedure. Second, we have generalized the method of argument used in establishing this improvement to provide a means by which any nondecresing sequence of constants can be rescaled so as to ensure control of the $FDP$. Finally, we have also described a procedure that controls the $FDR$, but only under an assumption on the joint distribution of the $p$-values.

In this article, we focused on the class of stepdown procedures. The alternative class of *stepup* procedures can be described as follows. Let

$$(5.1) \qquad \alpha_1 \leq \alpha_2 \leq \cdots \leq \alpha_s$$

be a nondecreasing sequence of constants. If $\hat{p}_{(s)} \leq \alpha_s$, then reject all null hypotheses; otherwise, reject hypotheses $H_{(1)}, \ldots, H_{(r)}$ where $r$ is the smallest index satisfying

$$(5.2) \qquad \hat{p}_{(s)} > \alpha_s, \ldots, \hat{p}_{(r+1)} > \alpha_{r+1}.$$

If, for all $r$, $\hat{p}_{(r)} > \alpha_r$, then reject no hypotheses. That is, a stepup procedure begins with the least significant $p$-value and continues accepting hypotheses as long as their corresponding $p$-values are large. If both a stepdown procedure and stepup procedure are based on the same set of constants $\alpha_i$, it is clear that the stepup procedure will reject at least as many hypotheses.

For example, the well-known stepup procedure based on $\alpha_i = i\alpha/s$ controls the $FDR$ at level $\alpha$, as shown by Benjamini and Hochberg [1] under the assumption that the $p$-values are mutually independent. Benjamini and Yekutieli [3] generalize their result to allow for certain types of dependence; also see Sarkar [14]. Benjamini and Yekutieli [3] also derive a procedure controlling the $FDR$ under no dependence assumptions. Romano and Shaikh [12] derive stepup procedures which control the $k$-$FWER$ and the $FDP$ under no dependence assumptions, and some comparisons with stepdown procedures are made as well.

## Acknowledgements

We wish to thank Juliet Shaffer for some helpful discussion and references.

## References

[1] Benjamini, Y. and Hochberg, Y. (1995). Controlling the false discovery rate: A practical and forceful approach to multiple testing. *J. Roy. Statist. Soc. Series B* **57**, 289–300. MR1325392

[2] Benjamini, Y. and Liu, W. (1999). A step-down multiple hypotheses testing procedure that controls the false discovery rate under independence. *J. Statist. Plann. Inference* **82**, 163–170. MR1736441

[3] Benjamini, Y. and Yekutieli, D. (2001). The control of the false discovery rate in multiple testing under dependency. *Ann. Statist.* **29**, 1165–1188. MR1869245

[4] Genovese, C. and Wasserman, L. (2004). A stochastic process approach to false discovery control. *Ann. Statist.* **32**, 1035–1061. MR2065197

[5] Hochberg, Y. and Tamhane, A. (1987). *Multiple Comparison Procedures.* Wiley, New York. MR0914493

[6] Holm, S. (1979). A simple sequentially rejective multiple test procedure. *Scand. J. Statist.* **6**, 65–70. MR0538597

[7] Hommel, G. (1983). Tests of the overall hypothesis for arbitrary dependence structures. *Biom. J.* **25**, 423–430. MR0735888

[8] HOMMEL, G. AND HOFFMAN, T. (1988). Controlled uncertainty. In *Multiple Hypothesis Testing* (P. Bauer, G. Hommel and E. Sonnemann, eds.). Springer, Heidelberg, 154–161.

[9] KORN, E., TROENDLE, J., MCSHANE, L. AND SIMON, R. (2004). Controlling the number of false discoveries: application to high-dimensional genomic data. *J. Statist. Plann. Inference* **124**, 379–398. MR2080371

[10] LEHMANN, E. L. AND ROMANO, J. (2005). Generalizations of the familywise error rate. *Ann. Statist.* **33**, 1138–1154. MR2195631

[11] PERONE PACIFICO, M., GENOVESE, C., VERDINELLI, I. AND WASSERMAN, L. (2004). False discovery rates for random fields. *J. Amer. Statist. Assoc.* **99**, 1002–1014. MR2109490

[12] ROMANO, J. AND SHAIKH, A. M. (2006). Stepup procedures for control of generalizations of the familywise error rate. *Ann. Statist.*, to appear. MR2195627

[13] SARKAR, S. (1998). Some probability inequalities for ordered $MTP_2$ random variables: a proof of Simes conjecture. *Ann. Statist.* **26**, 494–504. MR1626047

[14] SARKAR, S. (2002). Some results on false discovery rate in stepwise multiple testing procedures. *Ann. Statist.* **30**, 239–257. MR1892663

[15] SARKAR, S. AND CHANG, C. (1997). The Simes method for multiple hypothesis testing with positively dependent test statistics. *J. Amer. Statist. Assoc.* **92**, 1601–1608. MR1615269

[16] SIMES, R. (1986). An improved Bonferroni procedure for multiple tests of significance. *Biometrika* **73**, 751–754. MR0897872

[17] VAN DER LAAN, M., DUDOIT, S., AND POLLARD, K. (2004). Augmentation procedures for control of the generalized family-wise error rate and tail probabilities for the proportion of false positives. *Statist. Appl. Gen. Molec. Biol.* **3**, 1, Article 15. MR2101464